\documentclass{article}  

\usepackage{amsmath} 
\usepackage{amssymb}  
\usepackage{bm}

\makeatletter
\let\NAT@parse\undefined
\makeatother

\usepackage{amsthm}
\usepackage{graphicx}
\usepackage{enumerate}
\usepackage{mathtools}
\usepackage{xcolor} 

\usepackage{bbm}
\usepackage{algorithm}
\usepackage{algorithmic}
\usepackage{pgfplots}
\pgfplotsset{compat=newest}

\usepackage{tikz}
\usepackage{tkz-euclide}

\usepackage{enumerate}

\usepackage{mathrsfs}
\definecolor{ao(english)}{rgb}{0.0, 0.5, 0.0}
\usepackage[colorlinks, citecolor = {ao(english)}, linkcolor = {ao(english)}]{hyperref} 
\usepackage[capitalize,nameinlink]{cleveref}
\usepackage{siunitx}
\usepackage[sort,nocompress]{cite}
\usepackage{subcaption}
\usetikzlibrary{arrows,shapes,trees,calc,positioning,patterns,decorations.pathmorphing,decorations.markings,backgrounds}
\usetikzlibrary{matrix}

\usepgfplotslibrary{groupplots}

\newtheorem{thm}{Theorem}
\crefname{thm}{Theorem}{Theorems}

\newtheorem{prop}{Proposition}
\crefname{prop}{Proposition}{Propositions}

\newtheorem{lem}{Lemma}
\crefname{lem}{Lemma}{Lemmas}

\crefname{cor}{Corollary}{Corollaries}

\crefname{rem}{Remark}{Remark}

\crefname{ass}{Assumption}{Assumption}

\crefname{conj}{Conjecture}{Conjectures}

\crefname{defn}{Definition}{Definitions}

\crefname{prob}{Problem}{Problems}

\crefname{appl}{Application}{Applications}

\crefname{algorithm}{Algorithm}{Algorithms}
\crefname{paper}{Paper}{Papers}
\crefname{figure}{Figure}{Figures}
\crefname{section}{Section}{Sections}
\Crefname{section}{Section}{Sections}

\usepackage{dsfont}
\let\mathbb=\mathds

\newcommand{\Rmn}{\mathbb{R}^{n \times m}}

\newcommand{\conv}{\textnormal{conv}}
\newcommand{\rk}{\textnormal{rank}}
\newcommand{\svd}{\textnormal{svd}}
\newcommand{\diag}{\textnormal{diag}}

\newcommand{\argmin}{\operatornamewithlimits{argmin}}

\newcommand{\trace}{\textnormal{trace}}

\newcommand{\minmn}{\min \{ m,n \}}

\newcommand{\prox}{\textnormal{prox}}

\newcommand{\opts}{\star}

\newcommand{\card}{\textnormal{card}}
\newcommand{\dom}{\textnormal{dom}}

\newcommand{\normg}[1]{\|{}#1{} \|_{g}}

\newcommand{\normrgast}[1]{\| {} #1{} \|_{g,r\ast}}
\newcommand{\normrg}[1]{ \|{} #1{} \|_{g^D,r}}

\newcommand{\normA}[2]{\| {} #1{} \|_{ #2}}

\newcommand{\funcdom}{\Rmn \to \mathbb{R} \cup \lbrace \infty \rbrace}

\colorlet{FigColor1}{blue}
\colorlet{FigColor2}{red}
\colorlet{FigColor3}{ao(english)}
\colorlet{FigColor4}{orange}
\pgfplotsset{every axis plot/.append style={line width=1.5pt}}
\definecolor{bluebell}{rgb}{0.74, 0.83, 0.9}
\definecolor{airforceblue}{rgb}{0.36, 0.54, 0.66}

\crefformat{equation}{\textup{#2(#1)#3}}
\crefrangeformat{equation}{\textup{#3(#1)#4--#5(#2)#6}}
\crefmultiformat{equation}{\textup{#2(#1)#3}}{ and \textup{#2(#1)#3}}
{, \textup{#2(#1)#3}}{, and \textup{#2(#1)#3}}
\crefrangemultiformat{equation}{\textup{#3(#1)#4--#5(#2)#6}}%
{ and \textup{#3(#1)#4--#5(#2)#6}}{, \textup{#3(#1)#4--#5(#2)#6}}{, and \textup{#3(#1)#4--#5(#2)#6}}

\Crefformat{equation}{#2Equation~\textup{(#1)}#3}
\Crefrangeformat{equation}{Equations~\textup{#3(#1)#4--#5(#2)#6}}
\Crefmultiformat{equation}{Equations~\textup{#2(#1)#3}}{ and \textup{#2(#1)#3}}
{, \textup{#2(#1)#3}}{, and \textup{#2(#1)#3}}
\Crefrangemultiformat{equation}{Equations~\textup{#3(#1)#4--#5(#2)#6}}%
{ and \textup{#3(#1)#4--#5(#2)#6}}{, \textup{#3(#1)#4--#5(#2)#6}}{, and \textup{#3(#1)#4--#5(#2)#6}}

\crefdefaultlabelformat{#2\textup{#1}#3}

\renewcommand{\baselinestretch}{1}
\title{\LARGE Local Convergence of Proximal Splitting Methods\\for Rank Constrained Problems 
}

\author{Christian Grussler and Pontus Giselsson
	\thanks{The authors are members of the LCCC Linnaeus Center and the eLLIIT Excellence Center at Lund University. This work is financially supported by the Swedish Foundation for Strategic Research and the Swedish Research Council through the project 621-2012-5357. The authors are with the Department of Automatic Control, Lund University, Box 118, 22100 Lund, Sweden
		{\tt\small \{christiang,pontusg\}@control.lth.se}}%
	}

\begin{document}

	\maketitle
	\thispagestyle{empty}
	\pagestyle{empty}

	\begin{abstract}

	We analyze the local convergence of proximal splitting algorithms to solve optimization problems that are convex besides a rank constraint. For this, we show conditions under which the proximal operator of a function involving the rank constraint is locally identical to the proximal operator of its convex envelope, hence implying local convergence. The conditions imply that the non-convex algorithms locally converge to a solution whenever a convex relaxation involving the convex envelope can be expected to solve the non-convex problem.

	\end{abstract}
	
	\section{INTRODUCTION}
  Proximal splitting methods such as Douglas-Rachford splitting, the alternating direction method of multipliers, forward-backward splitting and many others (see~\cite{combettes2011proximal,bauschke2011convex,eckstein1992douglas,BoydDistributed,Gabay1976,Glowinski1975,douglas1956numerical}) are often used for solving large-scale convex optimization problems of the form
   \begin{equation}
   \label{eq:opt_org}
  \begin{aligned}
  & \underset{M}{\textnormal{minimize}}
  & & f_1(M)+ f_2(M),\\ 
  \end{aligned}
  \end{equation} 
  where $f_1$ and (or) $f_2$ have cheaply computable proximal mappings. 
  Since also many non-convex functions possess cheap proximal computations, there is a great interest in analyzing whether these iterates still converge to a solution. This paper focuses on analyzing the performance of splitting methods applied to problems, where $f_2$ is convex and 
   \begin{align}
  f_1(M) &:=  k (\|M\|) + \chi_{\rk(M) \leq r}(M), \label{eq:f_1}
  \end{align}
  is non-convex with 
  \begin{itemize}
  	\item $k(\cdot)$ being an increasing, convex function,
  	\item $\| \cdot \|$ being a unitarily invariant norm,
  	\item $\chi_{\rk(\cdot) \leq r}(\cdot)$ being the indicator function for matrices that have at most rank $r$.
  \end{itemize}
Analogously, one can consider vector-valued problems where the rank constraint is replaced by the cardinality constraint. Both problem types are very common within statistics, machine learning, automatic control and many more (see~\cite{grussler2016lowrank,grussler2016covariance,velu2013multivariate,vidal2016generalized,recht2010guaranteed,chandrasekaran2012convex,candes2006robust,hastie2015statistical}). 

Till this day, only special instances of solving this problem with proximal splitting methods have been analyzed \cite{hesse2014alternating,lewis1995convex,luke2013prox,hesse2013nonconvex,themelis2016forward}, mainly under the assumption that $f_2$ is the indicator function of an affine set and $k = 0$. In this paper, we deal with general convex functions $f_2$ and a large class of functions $f_1$, which allow us to provide an alternative analysis for showing local convergence. 

Letting $f_1^{\ast \ast}$ denote the bi-conjugate (convex envelope) of $f_1$, we show conditions under which the proximal operator to the non-convex function $f_1$ in \cref{eq:f_1} and its convex envelope $f_1^{**}$ (which was introduced in~\cite{grussler2016lowrank}) coincide. We translate these conditions to the setting of applying the Douglas-Rachford and forward-backward splitting algorithms to the non-convex problem \cref{eq:opt_org} with $f_1$ in \cref{eq:f_1}, and its optimal convex relaxation
  \begin{equation*}
\begin{aligned}
& \underset{M}{\textnormal{minimize}}
& & f_1^{\ast \ast}(M)+ f_2(M).
\end{aligned}
\end{equation*}
We show that the conditions imply local convergence of the non-convex splitting methods whenever all solutions to the convex relaxation are solutions to \cref{eq:opt_org}. Thus in many practical examples, there is no loss in directly using the non-convex algorithms. In fact, there are many examples where the non-convex methods can find a low-rank solution where the optimal convex relaxation fails. In other words, the non-convex algorithm can have low-rank limit points, whilst the convex has none, but not vice versa. This fact is  explicitly analysed for the case where $\|\cdot\|$ is the Frobenius norm and $k(\cdot) = (\cdot)^2$. 

Interestingly, we will see that unlike in the convex case, proximal splitting methods applied to \cref{eq:opt_org} and
   \begin{equation}
 \begin{aligned}
 & \underset{M}{\textnormal{minimize}}
 & & \gamma (f_1(M)+ f_2(M)),\\ 
 \end{aligned}
 \end{equation}
where $1 \neq \gamma > 0$, do not necessarily converge to the same limit points. Furthermore, the existence of a limit point as well as the the region of attraction in our local convergence result highly depend on the size of $\gamma$. On the one hand, if the optimal convex relaxation does not posses a low-rank solution, it is shown that $\gamma$ has to be chosen sufficiently small for a limit point to exists. On the other hand, in case of our guaranteed local convergence, the region of attraction grows with $\gamma$, i.e. for every initial point of the proximal algorithms there exists a sufficiently large $\gamma$ such that the algorithm converges.


Finally note that besides the ability of finding low-rank solutions when the convex relaxation fails, the non-convex algorithms are computationally more favourable, because the proximal computations of $f_1$ are significantly cheaper than those of the convex envelope $f_1^{\ast \ast}$ (see~\cite{grussler2017PhD}).

\section{Background}
The following notation for real matrices and vectors  $X=(x_{ij}) \in \Rmn$ is used in this paper. The non-increasingly ordered singular values of $X$, counted with multiplicity, are denoted by $$\sigma_1(X) \geq \dots \geq \sigma_{q}(X),$$ where $q:= \minmn$. Further, for $r \in \{1,\dots, q\}$ and $\sigma_r(X) \neq \sigma_{r+1}(X)$ we define the unique optimal rank-r approximation with respect to unitary invariant norms (see~\cite[Theorem~7.4.9.1]{horn2012matrix}) as
$$\svd_r(X) := \sum_{i=1}^r \sigma_i(X) u_i v_i^T,$$
where $X = \sum_{i=1}^{q} \sigma_i(X) u_i v_i^T$ is a singular value decomposition (SVD) of $X$. If $\sigma_r(X) = \sigma_{r+1}(X)$, then 
$$\svd_r(X) := \left \lbrace {\sum_{i=1}^r} \sigma_i(X) u_i v_i^T: X = \sum_{i=1}^{q} \sigma_i(X) u_i v_i^T \textnormal{ is an SVD of } X\right \rbrace.$$
Further, the inner-product for $X,Y \in \Rmn$ is defined by $$\langle X , Y \rangle := \sum_{i=1}^{m}\sum_{j=n}^{n} x_{ij} y_{ij} = \trace(X^TY).$$
\subsection{Norms}
A function $g: \mathbb{R}^q \to \mathbb{R}_{\geq 0}$ is called a \emph{symmetric gauge function} if
\begin{enumerate}[i.]
	\item $g$ is a~norm.
	\item $\forall x \in \mathbb{R}^{q}: g(|x|) = g(x)$, where $|x|$ denotes the element-wise absolute value.
	\item $g(Px) = g(x)$ for all permutation matrices $P\in\mathbb{R}^{q\times q}$ and all $x\in\mathbb{R}^q$.
\end{enumerate}
A norm $\normA{\cdot}{}$ on $\Rmn$ is \emph{unitarily invariant} if for all $X \in \Rmn$ and all unitary matrices $U$ and $V$ it holds that $\| U X V \| = \| X \|.$ Since all unitarily invariant norms on $\Rmn$ define a symmetric gauge function and vice versa (see~\cite{horn2012matrix}), 
we define $$\| \cdot \|_g := g(\sigma_1(\cdot),\dots,\sigma_{q} (\cdot)).$$ By~\cite{horn2012matrix} also the dual norm of $\normg{\cdot}$ is unitarily invariant and therefore it is associated with a symmetric gauge function $g^D$, i.e.
\begin{align*}
\| \cdot \|_{g^D} &:= \max_{\|X\|_g \leq 1} \langle \cdot , X \rangle= g^D(\sigma_1(\cdot),\dots,\sigma_q(\cdot)).
\end{align*}
For $r \in \{1,\ldots,q\}$, the \emph{truncated symmetric gauge functions} are given by $$g(\sigma_1,\dots,\sigma_r) := g(\sigma_1,\dots,\sigma_r,0,\dots,0).$$ Then, the so-called \emph{low-rank inducing norms} $\| \cdot \|_{g,r*}$ are defined in~\cite{grussler2016lowrank} as the dual norms of $$\| \cdot \|_{g^D,r} :=  g^D(\sigma_1(\cdot),\dots,\sigma_r(\cdot)).$$
The following properties have been shown in~\cite{grussler2016lowrank}.
\begin{lem} For all symmetric gauge functions $g:\mathbb{R}^q \to \mathbb{R}_{\geq 0}$ and $1 \leq r \leq q$ it holds that
	\label{lem:normrast}
\begin{align}
&\|M\|_g = \|M\|_{g,q\ast} \leq \dots \leq \|M\|_{g,1\ast}, \label{eq:norm_ineq}\\
&\rk(M)  \leq r \ \Rightarrow \ \|M \|_{g} = \|M \|_{g,r\ast}. \label{eq:rank_norm}
\end{align}
\end{lem} 
Finally, the Frobenius norm is given by $$\|X\|_{\ell_2}=\sqrt{\langle X , X \rangle} =  \sqrt{\sum_{i=1}^{q} \sigma_i^2(X)}.$$

\subsection{Functions}
The \emph{effective domain} of a function $f: \funcdom$ is defined as $$\dom f :=  \lbrace X \in \Rmn: f(X) < \infty \rbrace.$$ Then $f$ is said to be:
\begin{itemize}
	\item \emph{proper} if $\dom f \neq \emptyset$.
	\item \emph{closed} if for each $\alpha \in \mathbb{R}: \lbrace X \in \dom f: f(X) \leq \alpha \rbrace $ is a closed set.
\end{itemize}
A function $k:\mathbb{R} \to\mathbb{R}\cup \{ \infty \}$ is called \emph{increasing} if 
\begin{itemize}
	\item $x \leq y \ \Rightarrow \ k(x) \leq k(y) \text{ for all } x,y \in \mathbb{R},$
	\item $\exists \ x,y \in \mathbb{R}: k(x) < k(y)$.
\end{itemize}
The \emph{conjugate and bi-conjugate function} $f^\ast$ and $f^{\ast\ast}$ of $f$ are defined as $$f^\ast( \cdot ) :=   \sup_{X \in \Rmn} \left[ \langle X, \cdot \rangle - f(X) \right]$$ and $f^{\ast \ast} := (f^\ast)^{\ast}$. If $f: \mathbb{R} \to \mathbb{R} \cup \{ \infty \}$, then the \emph{monotone conjugate} is given by $$f^+(y) := \sup_{x \geq 0} \left[ \langle x, y \rangle - f(x) \right] \text{ for all } y\in\mathbb{R}.$$ 
The \emph{subdifferential} of $f$ in $X \in \dom f$ is defined as
\begin{align*}
\partial f(X) := \{G: f(Y) \geq f(X) + \langle G, Y-X\rangle \text{ for all } Y \},
\end{align*}
The proximal mapping of $f$ at $Z\in\Rmn$ is defined by 
\begin{align*}
\prox_{f} (Z) := \argmin_{M \in \Rmn}\left[f(M)+\frac{1}{2} \|M-Z\|_{\ell_2}^2 \right].
\end{align*}
Finally, for $\mathcal{S} \subset \Rmn$ the indicator function is defined as
\begin{align*}
\chi_{\mathcal{S}}(M) := \begin{cases}
0 & \text{if } M \in \mathcal{S},\\
\infty & \text{else}. 
\end{cases}
\end{align*}

\subsection{Optimal Convex Relaxation}
It is shown in~\cite{grussler2016lowrank} that the every low-rank inducing norm is the biconjugate (convex envelope) of \cref{eq:f_1} for different $\|\cdot\|_g$.
	\begin{prop}
	\label{prop:conj}
	Assume $k: \mathbb{R}_{\geq 0} \to \mathbb{R} \cup \{ \infty \}$ is an increasing closed  convex function, and let $f_{1} := k(\|\cdot\|_g) + \chi_{\rk(\cdot)  \leq r}$ be defined on $\Rmn$ with \linebreak $r \in \{1,\dots,\minmn \}$. Then,
	\begin{align}
	f^{\ast}_{1} &= k^{+}(\|\cdot\|_{g^D,r}), \\
	f^{\ast \ast}_{1} &= k(\|\cdot\|_{g,r\ast}).
	\end{align}
\end{prop}
These characterizations can be used to formulate Fenchel dual problems and \emph{optimal convex relaxations} to our rank constrained problems. This is shown in the following proposition, which is from~\cite{grussler2016lowrank}. 
	\begin{prop}
	\label{prop:opt_reg}
	Let $k: \mathbb{R}_{\geq 0} \to \mathbb{R} \cup \{ \infty \}$ and $f_2: \funcdom$ be proper, closed, convex functions with $r \in \{1,\dots,q\}$. Further let $k$ be increasing. Then,
			\begin{align}
			\inf_{\stackrel{M \in \Rmn}{\rk(M) \leq r}} \left[k(\|M\|_g) + f_2(M) \right]& \label{eq:prim} \\ 
			&\geq  -\min_{D \in \Rmn} \left[k^{+}(\|D\|_{g^D,r})  + f_2^\ast(-D) \right] \label{eq:dual}\\ 
			&= \min_{M \in \Rmn} \left[k(\|M\|_{g,r\ast}) + f_2(M) \right] \label{eq:bidual}.
			\end{align}
	If $M^\opts$ solves \cref{eq:bidual} such that $\rk(M^\opts) \leq r$, then equality holds, and $M^\opts$ is also a solution to \cref{eq:prim}.
\end{prop}
	
\section{Theoretical Results}
In this section we derive the theoretical results that are needed for our convergence analysis  in~\cref{sec:conv}. The proofs to these results are given in the appendix.
\begin{thm}
	\label{thm:prox}
Let $Z \in \Rmn$, $f_1 := k (\normg{\cdot}) + \chi_{\rk(\cdot) \leq r}$ and $\gamma > 0$, where $k: \mathbb{R}_{\geq 0} \to \mathbb{R} \cup \{ \infty \}$ is a proper, closed and increasing convex function and $r \in \{1,\dots,q\}.$ Then for all $P \in \svd_r(Z)$ it holds that
\begin{align*}
\prox_{\gamma k (\normg{\cdot})}(P) &= \prox_{\gamma k (\normrgast{\cdot})}(P) \in \prox_{\gamma f_1}(Z).
\end{align*}
Moreover, let $$M^c := \prox_{\gamma k (\normrgast{\cdot})}(Z),$$
%
then the following are equivalent:
\begin{enumerate}[i.]
	\item $M^c = \prox_{\gamma f_1}(Z),$ \label{item:prox_eq}
	\item $\rk(M^c) \leq r,$ \label{item:prox_rank}
	\item $\sigma_{r+j}(Z-M^c) = \sigma_{r+j}(Z)$ for all $j \in \{1,\dots,q-r\}$.\label{item:prox_sig_eq}
	\item $\sigma_r(Z-M^c) \geq \sigma_{r+1}(Z).$ \label{item:prox_sig}
\end{enumerate}
\end{thm}
Computing the prox of the non-convex function $f_1$ at $Z$, reduces to evaluating the convex prox of either $k(\|\cdot\|)_g$ or the convex envelope $f_1^{**}=k(\|\cdot\|_{g,r*})$ at $P\in\svd_r(Z)$. Therefore, only the first $r$ singular values and vectors are needed to compute the non-convex prox. This can be compared to the prox of the convex envelope $f_1^{**}$ at $Z$, where all singular values and vectors might be needed. To compute the prox of $k(\|\cdot\|)_g$ is cheaper than computing the prox of $f_1^{**}=k(\|\cdot\|_{g,r*})$, except for rank-$r$ matrices, see~\cite{grussler2016lowrank}. Therefore it is often much cheaper to evaluate the prox of the non-convex function $f_1$ than of its convex envelope $f_1^{**}=k(\|\cdot\|_{g,r*})$. 

In order to relate \cref{thm:prox} to the solutions of \cref{eq:prim} and \cref{eq:bidual}, the following results, which are proven in \cref{proof:rank,proof:rankDM} respectively, will be needed. 
\begin{lem}
	\label{lem:rank}
	Let $D \in \Rmn$ and $1 \leq r \leq q$. Assume that $$\sigma_r(D) = \dots = \sigma_{r+s}(D),$$ where either $s = q-r$ or $\sigma_{r+s}(D) \neq \sigma_{r+s+1}(D)$ for some $s\geq0$. Then all $M \in \partial \normrg{D}$ fulfill that 
	$$\rk(M) \leq r+s.$$
	Moreover, if $\sigma_r(D) = 0$, then $\rk(M) \leq r$. 
\end{lem}
\begin{prop}
	\label{prop:rankDM}
	Let $D^\opts, M^\opts \in \Rmn $ be solutions to \cref{eq:dual} and \cref{eq:bidual}, respectively. Assume $$\sigma_{r}(D^\opts) = \dots = \sigma_{r+s}(D^\opts) \neq 0,$$ where either $s = q-r$ or $\sigma_{r+s}(D^\opts) \neq \sigma_{r+s+1}(D^\opts)$ for some $s \geq 0$. Further, if $\sigma_{r}(D^\opts) = 0$ let $s =0$.
	Then, $$\rk(M^\opts) \leq r+s.$$
	In particular, if there exists a solution $D^\opts$ to \cref{eq:dual} such that $\sigma_{r}(D^\opts) \neq \sigma_{r+1}(D^\opts)$ or $\sigma_r(D^\opts) = 0$, then all solutions to \cref{eq:bidual} are solutions to \cref{eq:prim}.
\end{prop}

\section{Convergence Analysis}
\label{sec:conv}
Next it is discussed how \cref{thm:prox,prop:rankDM} can be used to show local convergence of proximal splitting algorithms applied to problems of the form
\begin{equation}
\label{eq:opt_conv}
\begin{aligned}
& \underset{M}{\textnormal{minimize}}
& &  k (\normg{M}) + \chi_{\rk(M) \leq r}(M) + f_2(M),\\ 
\end{aligned}
\end{equation} 
where $f_2$ is a convex function with cheaply computable proximal mapping and $k$ an convex, increasing function. To illustrate and support our analysis, let us first recap the following two well-known proximal splitting algorithms applied to~\cref{eq:opt_org}.

\subsection*{Douglas-Rachford Splitting}
The Douglas-Rachford splitting method is one of the most well-known splitting algorithms for solving large-scale convex problems~\cite{douglas1956numerical,lions1979splitting,eckstein1992douglas,combettes2011proximal}. In fact, the well-known alternating direction methods of multipliers (ADMM) is a special case of this algorithm (see~\cite{Glowinski1975,Gabay1976,BoydDistributed}). The Douglas-Rachford iterations are given by
\begin{subequations}
	\begin{align}
	X_{k} &= \prox_{\gamma f_1}(Z_{k-1}), \label{eq:DR_X}\\
	Y_{k} &= \prox_{\gamma f_2}(2X_k-Z_{k-1}), \label{eq:DR_Y}\\
	Z_{k} &= Z^{k-1}+\rho(Y_k-X_k),\label{eq:DR_Z}
	\end{align}
\end{subequations}
where $\gamma > 0$ and $0 < \rho < 2$. For convex $f_1$ and $f_2$, $X^k$ and $Y^k$ converge towards an identical solution of \cref{eq:opt_org} and $\{\|Z_k - Z^\opts\|_{\ell_2}\}_{k \in \mathbb{N}}$ is non-increasing, where   $Z^\opts := \lim_{k\to\infty} Z_k$. (see~\cite{douglas1956numerical,lions1979splitting,eckstein1992douglas}).

\subsection*{Forward-Backward Splitting}
Another popular splitting methods is the so-called forward-backward splitting algorithm (see~\cite{combettes2011proximal,bauschke2011convex,levitin1966constrained,rockafellar1970convex}). In this case, $f_2$ is assumed to be differentiable with Lipschitz continuous gradient, i.e. for all $X, Y \in \Rmn$
\begin{equation*}
\|\nabla f_2 (Y) - \nabla f_2(X)\|_{\ell_2} \leq L \|Y-X\|_{\ell_2}. 
\end{equation*}
Then the forward-backward iterations are given by
\begin{align*}
Z_{k} &= X_k - \gamma \nabla f_2(X_{k-1}),\\
X_{k} &= \prox_{\gamma f_1}(Z_k),
\end{align*}
where $ 0 < \gamma < \frac{2}{L}$. Also here if $f_1$ and $f_2$ are convex, then it can be shown that $X^k$ converges towards a solution of \cref{eq:opt_org} and 
$\{\|Z_k - Z^\opts\|_{\ell_2}\}_{k \in \mathbb{N}}$ is non-increasing  
with  $Z^\opts := \lim_{k\to\infty} Z_k$.

\subsection*{Local Convergence}
One of the steps in the above two methods (and many other operator splitting methods) when applied to solve \cref{eq:opt_org} is
\begin{align*}
X_k = \prox_{\gamma f_1}(Z_{k-1}).
\end{align*}
If $f_1$ and $f_2$ are convex, then $X_k$ converges to a solution of \cref{eq:opt_org} in both methods and $\{\|Z_k - Z^\opts\|_{\ell_2}\}_{k \in \mathbb{N}}$ is a non-increasing sequence, where  $Z^\opts := \lim_{k\to \infty} Z_k$. Next, we will show that the latter and \cref{thm:prox} imply local convergence of proximal splitting algorithms applied to the non-convex problem in~\cref{eq:opt_conv}. 

In the following we will refer to a proximal splitting algorithm applied to the optimal convex relaxation in \cref{eq:bidual}, which is restated here,
\begin{equation}
\label{eq:opt_convex}
\begin{aligned}
& \underset{M}{\textnormal{minimize}}
& &  k (\normrgast{M}) + f_2(M),\\ 
\end{aligned}
\end{equation} 
as the \emph{convex splitting algorithm} with iterates 
\begin{align*}
M^c_k = \prox_{\gamma k (\normrgast{\cdot})}(Z_k).
\end{align*}
Correspondingly, if the algorithm is applied to \cref{eq:opt_conv}, i.e.  $f_1 = k (\|\cdot\|_g) + \chi_{\rk(\cdot) \leq r}$ , we speak of the \emph{non-convex splitting algorithm} with iterates
\begin{align*}
M^n_k = \prox_{k (\|\cdot\|_g) + \chi_{\rk(\cdot) \leq r}}(Z_k).
\end{align*}
Let us assume that $M^\opts$ is a solution to~\cref{eq:opt_convex} with
\begin{enumerate}[i.]
	\item $Z^\opts = \lim_{k \to \infty } Z_k,$ 
	\item $M^\opts = \lim_{k \to \infty } M^c_k,$
	\item $\sigma_r(Z^\opts-M^\opts) > \sigma_{r+1}(Z^\opts).$
\end{enumerate} 
By (firm) nonexpansiveness of $\prox_{\gamma k (\normrgast{\cdot})}(Z)$ 
and the continuity of the singular values (see~\cite[Corollary~4.9]{stewart1990matrix}), \cref{thm:prox} implies that\begin{align*}
M^c_0= M^n_0
\end{align*}
for all $Z_0 \in B_{\varepsilon}(Z^\opts) := \{Z:\|X-Z^\opts\|_{\ell_2} < \varepsilon \}$, where  $\varepsilon:= \sigma_r(Z^\opts-M^\opts) - \sigma_{r+1}(Z^\opts) > 0$. Thus, since   $\{\|Z_k^c - Z^\opts\|_{\ell_2}\}_{k \in \mathbb{N}}$ is non-increasing, it follows that
\begin{align*}
\forall k \geq 0, \ Z_0 \in B_{\varepsilon}(Z^\opts): \ M^c_k= M^n_k.
\end{align*}
This proves the \emph{local convergence} of the non-convex algorithm if
$\sigma_r(Z^\opts-M^\opts) > \sigma_{r+1}(Z^\opts)$. 

We will conclude this section by linking this condition to the solution set of~\cref{eq:opt_convex}, which is the same as \cref{eq:bidual}. A necessary optimality condition for solving \cref{eq:dual} and \cref{eq:bidual} is that (see~\cite[Theorem~7.12.1]{luenberger1968optimization} and \cite[Theorem~23.5.]{rockafellar1970convex})
\begin{align*}
D^\opts\in\partial_M \left. k(\|M\|_{g,r*}) \right|_{M = M^\opts}.
\end{align*}
Now, relating this to the optimality condition of the convex prox computation:
\begin{align*}
0\in\partial_M \left. k(\|M\|_{g,r*}) \right|_{M = M^\opts}+\gamma^{-1}(M^\opts-Z^\opts),
\end{align*}
implies that
\begin{align*}
D^\opts = \gamma^{-1}(Z^\opts - M^\opts)
\end{align*}
is a solution to the dual problem \cref{eq:dual}, i.e.,
\begin{equation}
\label{eq:opt_dual}
\begin{aligned}
& \underset{D}{\textnormal{minimize}}
& &  k^{+}(\|D\|_{g^D,r})  + f_2^\ast(-D).\\ 
\end{aligned}
\end{equation}
By \cref{thm:prox} we can conclude that $\sigma_{r+1}(Z^\opts) = \gamma \sigma_{r+1}(D^\opts)$ if $\sigma_r(Z^\opts-M^\opts) \geq \sigma_{r+1}(Z^\opts)$. Hence, \cref{prop:rankDM} implies the local convergence of non-convex proximal splitting algorithms, if there exists a solution $D^\opts$ to \cref{eq:opt_dual} such that
\begin{align}
\sigma_r(D^\opts) \neq \sigma_{r+1}(D^\opts) \ \text{or} \ \sigma_r(D) 
\label{eq:cond_lowrank}
 = 0.\end{align}
This condition insures, by \cref{prop:rankDM}, that \cref{eq:opt_convex} has only solutions of at most rank $r$. Note that if \cref{eq:opt_convex} has solutions of rank larger than $r$, then a convex algorithm cannot be expected to find solutions of rank $r$, despite their possible existence. This is because the solution set of a convex problem is a convex set.

In other words, non-convex proximal splitting methods locally converge to a solution of \cref{eq:opt_conv}, whenever one can expect to find such a solution by solving \cref{eq:opt_convex}. Moreover, the \emph{region of attraction} to $Z^\opts$ contains the ball $B_{\varepsilon}(Z^\opts)$ with $$\varepsilon = \gamma ( \sigma_{r}(D^\opts) - \sigma_{r+1}(D^\opts)) = \sigma_r(Z^\opts-M^\opts) - \sigma_{r+1}(Z^\opts) > 0.$$
This means that for each initial point $Z_0$ there exists a $\gamma >0$ that guarantees the convergence to $Z^\opts$. Finally, numerical experiments indicate that the non-convex algorithms can also find rank-r solutions to \cref{eq:opt_convex,} despite the fact that \cref{eq:opt_convex} may have higher rank solutions. 

\section{Douglas-Rachford Limit Points}
In the following, let us compare the Douglas-Rachford limit points to the optimal convex relaxation (convex Douglas-Rachford) with the limit points of the non-convex Douglas-Rachford for problems \cref{eq:opt_org} where
   \begin{align*}
f_1(M) &:=  \frac{1}{2}\|\cdot\|_{\ell_2}^2 + \chi_{\rk(M) \leq r}(M).
\end{align*}
Using completion of squares and the well-known Schmidt-Mirsky Theorem (see~\cite[Theorem~7.4.9.1]{horn2012matrix}), we get that
\begin{align}
\prox_{\gamma f_1}(Z) &= \argmin_{\stackrel{M \in \Rmn}{\rk(M)\leq r}}\left( \frac{\gamma}{2} \|M\|_F^2 + \frac{1}{2} \|M-Z\|_F^2 \right) \notag \\ 
&= \argmin_{\stackrel{M \in \Rmn}{\rk(M)\leq r}}\left(\frac{\gamma + 1}{2} \|M\|_F^2 - \langle Z,M \rangle \right) \notag \\
&= \argmin_{\stackrel{M \in \Rmn}{\rk(M)\leq r}} \left \| \frac{Z}{\gamma +1}-M \right\|_F^2 \notag \\
&= \frac{1}{1+\gamma}\svd_r \left({Z} \right). \label{eq:prox_nonDR}
\end{align}
This allows us to derive the following comparative result on the limit points of the convex and non-convex Douglas-Rachford, which is proven in \cref{proof:fixpoint}. 
\begin{thm}\label{thm:fixpoint}
	Let $X^\opts \in \Rmn$ with $\rk(X^\opts) \leq r$ and $\gamma > 0$. Then $X^\opts$ is a limit point of the \emph{convex (non-convex) Douglas-Rachford splitting} iterate \cref{eq:DR_X} if and only if there exists $R \in \Rmn$ such that
	\begin{equation*}
	R^T X^\opts = 0, \quad X^\opts R^T = 0, \quad -X^\opts - R \in \partial g(X^\opts),
	\end{equation*}
	and in the 
	\begin{itemize}
		\item convex case: ${\sigma_1(R) \leq \sigma_r(X^\opts)},$
		\item non-convex case: ${\sigma_1(R) \leq (1+\gamma^{-1})\sigma_r(X^\opts)}.$
	\end{itemize}
\end{thm}
\cref{thm:fixpoint} verifies what has been discussed in the end of previous section that all limit points of the convex Douglas-Rachford are limit points to the non-convex Douglas-Rachford, but not vice versa. More importantly, it shows the importance of choosing a feasible $\gamma$. In the presence of a duality gap in~\cref{eq:dual}, \cref{thm:fixpoint} implies that if $\gamma$ is chosen too large, then the non-convex Douglas-Rachford may not posses a limit point, but choosing $\gamma$ sufficiently small can help to gain convergence. Analytical examples where this applies have been studied in~\cite{grussler2017PhD} and a numerical example is given in the next section. This is very much in contrast to the convex case, where convergence is independent of $\gamma$. Finally note that by choosing $\gamma$ just small enough for a limit point to exist, the problem of multiple limit points may be avoided and thus making the algorithm independent of the initialization. Similar derivations can be carried out for all $f_1$ in the form of \cref{eq:f_1}. 

\pagebreak
\section{Example}
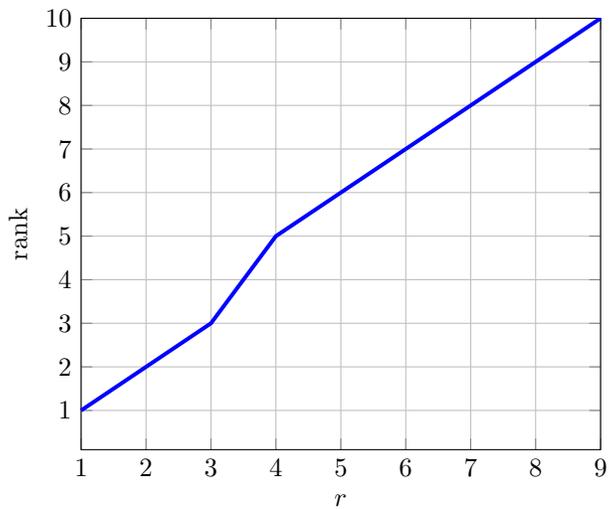
\begin{figure}
	\centering
	\def\factor{.7}%
	\begin{tikzpicture}
\begin{axis}[xlabel=$r$,
	ylabel=$\rk$,
	xmin = 1,
	xmax = 9,
	grid = both,
	width = \factor \textwidth,
	ymax = 10,
	xtick = {1,...,9},
	ytick = {1,...,10},
	scaled ticks=false, 
	]
	
	\addplot+[color = FigColor1,mark = none] file{rank_conv.txt};
	\label{hankel_pres_error:lower}

	\end{axis}

	\end{tikzpicture}
	\caption{Hankel matrix approximation -- Rank of the solutions to the optimal convex relaxation. \label{fig:rank}}
\end{figure}

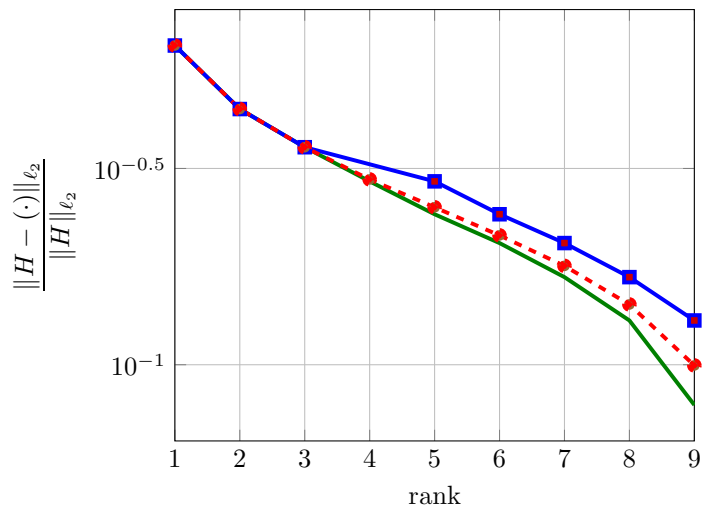
\begin{figure}
	\centering
	\def\factor{.7}%
	\begin{tikzpicture}
	\begin{axis}[xlabel=$\rk$,
	ylabel=$\dfrac{\|H-(\cdot)\|_{\ell_2}}{\|H\|_{\ell_2}}$,
	xmin = 1,
	xmax = 9,
	grid = both,
	width = \factor \textwidth,
	ymode =  log,
	xtick = {1,...,9},
	scaled ticks=false, 
	]

	\addplot+[color = FigColor3,mark = none] file{lower.txt};
	\label{line:lower}
	
	\addplot+[color = FigColor1] file{err_conv.txt};
	\label{line:conv}
	\addplot+[color = FigColor2, dashed] file{err_ndr.txt};
	\label{line:ndr}
	\end{axis}

	\end{tikzpicture}
	\caption{Hankel matrix approximation -- Relative errors of the approximations obtained by \ref{line:conv} the optimal convex relaxation and \ref{line:ndr} non-convex Douglas-Rachford. \ref{line:lower} indicates the lower bound obtained by the optimal convex relaxation. \label{fig:err}}
\end{figure}

Within many areas such as automatic control, the rank of a Hankel operator/matrix is crucial, because it determines the order of a linear dynamical system. Whereas, the celebrated Adamyan-Arov-Krein theorem (see~\cite{antoulas2005approximation}) answers the question of optimal low-rank approximation of infinite dimensional Hankel operators, the following finite dimensional case is still unsolved:
\begin{equation*}
	\begin{aligned}
	& \underset{M}{\textnormal{minimize}}
	& & \|H-M\|_{\ell_2}^2 \\
	& \textnormal{subject to}
	& & \rk(M) \leq r,\\ 
	& & & M \in \mathcal{H},
	\end{aligned}
\end{equation*}	
where $H \in \mathcal{H} := \lbrace X \in \mathbb{R}^{n \times n} : X \textnormal{ is Hankel} \rbrace$. In the following, we show how non-convex Douglas-Rachford splitting performs on this problem class in comparison with the optimal convex relaxation. To this end, we rewrite the problem in the view of \cref{eq:opt_convex} and \cref{eq:bidual} as 
\begin{equation*}
\begin{aligned}
& \underset{M}{\textnormal{minimize}}
& &  \|M\|_{\ell_2}^2 +  \chi_{\rk(M) \leq r}(M) + f_2(M),
\end{aligned}
\end{equation*}	
where $f_2(M) := 2 \langle M,H \rangle + \|H\|_{\ell_2}^2 + \chi_{\mathcal{H}}(M)$. For our numerical experiments we use
\begin{equation*}
H := \begin{tikzpicture}[baseline=(current bounding box.center)]
\matrix (m) [matrix of math nodes,
nodes in empty cells,
right delimiter={)},
left delimiter={(}]{
	1  	& 	1 	& 	  	&   	& 	1 	& 	1  \\
	1  	& 		& 		& 		&  	    & 	0	\\
	&		&		&		&		&		\\
	&		&		&		&		& 		\\
	1	&		&		&		&		& 	0  	\\
	1   &  	0   &		&		&	0	&	0	\\
} ;
\newdimen\L
\L = .8 pt
\draw[loosely dotted, line width = \L] (m-1-2)-- (m-1-5);
\draw[loosely dotted, line width = \L] (m-6-2)-- (m-6-5);
%
\draw[loosely dotted, line width = \L] (m-2-1)-- (m-5-1);
\draw[loosely dotted, line width = \L] (m-2-6)-- (m-5-6);
%
\draw[loosely dotted, line width = \L] (m-5-1)-- (m-1-5);
\draw[loosely dotted, line width = \L] (m-6-1)-- (m-1-6);
\draw[loosely dotted, line width = \L] (m-6-2)-- (m-2-6);


\end{tikzpicture} \in \mathbb{R}^{10 \times 10}.
\end{equation*}
The non-convex Douglas-Rachford uses $\gamma =1$ and is initialized with $Z_0 = 0$ for all $r \in \{1,\dots,9\}$. The ranks of the solutions to the optimal convex relaxation are shown in \cref{fig:rank}. We observe that only for $r = \{1,2,3\}$ the convex relaxation manages to find guaranteed solutions to the non-convex problem. In contrast, the non-convex Douglas Rachford converges for all $r$. \cref{fig:err} shows the relative errors of these solutions and the (sub-optimal) solutions to the convex relaxation as well as the lower bound that is provided by the convex relaxation (see~\cref{prop:opt_reg}). Note that the convex relaxation is not able to obtain a sub-optimal solution of rank $4$. From~\cref{fig:err} it can be seen that the non-convex solutions for $r = \{1,2,3\}$ coincide with the convex solutions, just as our local convergence guarantee suggests. However, for all other $r$, the non-convex approximations outperform the sub-optimal solutions of the convex relaxation. Finally, is has been observed that, if one chooses $\gamma$ sufficiently large, the non-convex Douglas-Rachford does not converge for $r > 3$. This can be explained through \cref{thm:fixpoint}. 

\section{Conclusion}
We have shown conditions under which the proximal mapping of the non-convex function \cref{eq:f_1} coincides with the proximal mapping of its convex envelope. This allowed us to state conditions under which the non-convex and convex Douglas-Rachford methods and forward-backward methods coincide. This, in turn, guarantees local convergence of the non-convex methods in these situations. Furthermore, we have provided a comparison between the convex and non-convex Douglas-Rachford limit points for common instance of the squared Frobenius norm. Unlike in the convex case, this has demonstrated that scaling the problem may have significant impact. Finally, we discussed a numerical example in which a non-convex method converges also when the stated assumptions do not hold. In those situations, the quality of the solution from the non-convex algorithm was better than the solution obtained by the optimal convex relaxation.

\bibliography{refopt,refpos}
\appendix
\section{Appendix}
	\subsection{Proof to \cref{thm:prox}}
	\begin{proof}
		For $M \in \Rmn$ and $1 \leq j \leq q$, let us define $$\Sigma_j(M) := \diag(\sigma_1(M),\dots,\sigma_j(M),0,\dots,0) \in \mathbb{R}^{q \times q}.$$ By \cite[Corollary~2.5.]{lewis1995convex} and the unitary invariance of $k(\normg{\cdot})$, it can be seen that $X := \prox_{\gamma k (\normg{\cdot})}(P)$ and $P$ have simultaneous SVDs, i.e. if $P = \sum_{i=1}^r \sigma_i(P)u_iv_i^T$, then $\prox_{\gamma k (\normg{\cdot})}(P) = \sum_{i=1}^r \sigma_i(X)u_iv_i^T$. Hence,	
		\begin{align*}
		&\argmin_{{M\in \Rmn}}\left[\gamma k (\normg{M})+\frac{1}{2} \|M-P\|_{\ell_2}^2 \right]\\
		&=\argmin_{M\in \Rmn}\left[\gamma k (\normg{\Sigma_q(M)})+\frac{1}{2} \|\Sigma_q(M)-\Sigma_r(P)\|_{\ell_2}^2 \right].
		\end{align*}
		Further, \cite[Theorem~7.4.8.4.]{horn2012matrix} implies that
		\begin{align*}
		\normg{\Sigma_q(M)} &\geq \normg{\Sigma_r(M)},\\
		\|\Sigma_q(M)-\Sigma_r(Z)\|_{\ell_2} &\geq \|\Sigma_r(M)-\Sigma_r(Z)\|_{\ell_2},
		\end{align*}
		for all $M \in \Rmn$, which yields that
		\begin{align*}
		\argmin_{M\in \Rmn}&\left[\gamma k (\normg{\Sigma_q(M)})+\frac{1}{2} \|\Sigma_q(M)-\Sigma_r(P)\|_{\ell_2}^2 \right]\\
		&= \argmin_{M\in \Rmn}\left[\gamma k (\normg{\Sigma_r(M)})+\frac{1}{2} \|\Sigma_r(M)-\Sigma_r(P)\|_{\ell_2}^2 \right]\\
		&= \argmin_{M\in \Rmn}\left[\gamma k (\normg{\Sigma_r(M)})+\frac{1}{2} \|\Sigma_r(M)-\Sigma_q(Z)\|_{\ell_2}^2 \right]\\
		&=\argmin_{\stackrel{M\in \Rmn}{\rk(M) \leq r}}\left[\gamma k (\normg{\Sigma_q(M)})+\frac{1}{2} \|\Sigma_q(M)-\Sigma_q(Z)\|_{\ell_2}^2 \right]\\
		&\in \argmin_{\stackrel{M\in \Rmn}{\rk(M) \leq r}}\left[\gamma k (\normg{M})+\frac{1}{2} \|M-Z\|_{\ell_2}^2 \right],
		\end{align*}
		where the last equality and the inclusion follow by \cite[Corollary 7.4.1.3.]{horn2012matrix}, \cite[Corollary~2.5.]{lewis1995convex} and the unitary invariance of $f_1$. This proves that $$\prox_{\gamma k (\normg{\cdot})}(P) \in \prox_{\gamma f_1}(Z).$$ Moreover, by \cref{eq:norm_ineq} it follows that $\rk(\prox_{\gamma k (\normg{\cdot})}(P)) \leq r$ implies $$\prox_{\gamma k (\normg{\cdot})}(P) = \prox_{\gamma k (\normrgast{\cdot})}(P).$$  
		By the extend Moreau decomposition (see~e.g.~\cite{BoydDistributed}) and \cref{prop:conj} it holds that 
		\begin{align*}
		M^c = Z - Y^c,
		\end{align*}
		where $Y^c := \gamma \prox_{\gamma^{-1} k^{+} (\normrg{\cdot})}(\gamma^{-1}Z)$.
		As before, $M^c$, $Z$ and $Y^c$ can be shown to have simultaneous SVDs which is why
		\begin{align}
		\Sigma_q(M^c) = \Sigma_q(Z) - \Sigma_q(Y^c).
		\end{align}
		Thus $\rk(M^c) \leq r$ if and only if 	
		$\sigma_{j}(Y^c) = \sigma_{j}(Z)$ for $r+1 \leq j \leq q$. Since, $\normrg{\cdot}$ only depends on $\sigma_1(Y^c),\dots,\sigma_r(Y^c)$, this is equivalent to $$\sigma_{r}(Y^c) \geq \sigma_{r+1}(Z).$$ This shows the equivalence between \cref{item:prox_rank,item:prox_sig_eq,item:prox_sig}. Finally note that this is also equivalent to $$\Sigma_q(M^c) = \prox_{\gamma k (\normrgast{\cdot})}(\Sigma_r(Z)).$$  Since $M^c$ is unique, this can only be true if  $\sigma_r(Z) \neq \sigma_{r+1}(Z)$ and thus $M^c = \prox_{\gamma f_1}(Z)$, which concludes the proof.	
	\end{proof}
	\subsection{Proof to \cref{lem:rank}}
        \label{proof:rank}
	\begin{proof}
		Let $D = \sum_{i=1}^{q} \sigma_i(D)u_iv_i^T$ be an SVD of $D$ and $\bm{\sigma}(D) \in \mathbb{R}^q$ the corresponding vector of singular values. Further, let 
		for all $\mathcal{I} \subset \{1,\dots,q\}$, \linebreak $\Pi_{\mathcal{I}}: \mathbb{R}^q \to \mathbb{R}^q$ be defined as  
		\begin{align*}
		\Pi_{\mathcal{I}}(x)_i &:= \begin{cases}
		x_i & \text{if } i \in \mathcal{I}\\
		0   & \text{if } i \notin \mathcal{I}.
		\end{cases}\\
		\end{align*}
		By \cite[Theorem~2]{watson1992characterization} it holds that
		\begin{align}
		\partial \normrg{D} = \left\lbrace \sum_{i=1}^{q} d_i u_i v_i^T: d \in \partial \normrg{\bm{\sigma}(D)}\right\rbrace. 
		\end{align}
		Next we show that $d_{r+s+1} = \dots = d_{q} =0$. Letting $\card(\cdot)$ denote the cardinality, it follows from \cite[Theorem~7.4.8.4.]{horn2012matrix} that
		\begin{align*}
		\normrg{x} = \max_{\card(\mathcal{I}) = r} g^D(\Pi_{\mathcal{I}}(x)).
		\end{align*}
		and therefore by \cite[Corollary~VI.4.3.2]{hiriart2013convex}
		\begin{align*}
		\partial \normrg{\bm{\sigma}(D)}= \conv (\lbrace \partial g^D(\Pi_{\mathcal{I}}(\bm{\sigma}(D))): \card(\mathcal{I}) = r, \\
		g^D(\Pi_{\mathcal{I}}(\bm{\sigma}(D))) = \normrg{\bm{\sigma}(D)} \rbrace),
		\end{align*}
		where $\conv(\cdot)$ denotes the convex hull. However, \cite[Theorem~7.4.8.4.]{horn2012matrix} implies that $\mathcal{I} \subset \{1,\dots,r+s\}$ if $$g^D(\Pi_{\mathcal{I}}(\bm{\sigma}(D))) = \normrg{\bm{\sigma}(D)}.$$ In this case, $g^D(\Pi_{\mathcal{I}}(x))$ only depends on variables $x_1,\dots,x_{r+s}$, which is why for all $d \in \partial g^D(\Pi_{\mathcal{I}}(\bm{\sigma}(D)))$ and hence all $d \in \partial \normrg{\bm{\sigma}(D)}$ it holds that $$d_{r+s+1} = \dots = d_{q}.$$ This proves the first claim. Then the second claim follows by the continuity of the subdifferentials (see~\cite[Theorem~24.4.]{rockafellar1970convex}).
			\end{proof}
		\subsection{Proof to \cref{prop:rankDM}}
                \label{proof:rankDM}
		\begin{proof}
			By Fenchel's Duality Theorem (see~\cite[Theorem~7.12.1]{luenberger1968optimization}) and \cite[Theorem~23.5.]{rockafellar1970convex} it follows that
			\begin{align*}
			M^\opts \in \left. \partial_D k^{+}(\|D\|_{g^D,r}) \right|_{D = D^\opts}.
			\end{align*}
			Since $k^+$ is increasing (see~[p.~111]\cite{rockafellar1970convex}) 
			it holds by~\cite[Theorem~VI.4.3.2]{hiriart2013convex} that there exist
			\begin{align*}
			p \in \partial k^{+}(\|D^\opts\|_{g^D,r}) \quad \textnormal{and} \quad N \in \partial \|D^\opts\|_{g^D,r}
			\end{align*}
			such that $ M^\opts = p N$. Then invoking \cref{lem:rank} proves the claims.
		\end{proof}
		\subsection{Proof to \cref{thm:fixpoint}}
		\label{proof:fixpoint}
\begin{proof}
	In the convex case, $X^\opts$ is a limit point of \cref{eq:DR_X} if and only if
	$X^\opts$ solves \cref{eq:bidual}. Letting $D^\opts$ be a solution to \cref{eq:dual}, this is equivalent to
	\begin{align*}
	X^\opts \in \svd_r(D^\opts), \quad -D^\opts \in \partial g(X^\opts),
	\end{align*}
	by \cite[Theorem~3]{grussler2016low} and \cite[Theorems~23.5 and 31.1]{rockafellar1970convex}. Thus, defining $R := D^\opts - X^\opts$ proves the equivalence for the convex case. 
	
	Next let $X^\opts$ be a limit point to the non-convex Douglas Rachford splitting algorithm. By \cref{eq:DR_X}, \cref{eq:DR_Y}, \cref{eq:prox_nonDR} and \cite[Theorems~23.5 and 27.1]{rockafellar1970convex}, this holds if and only if there exists a $Z^\opts \in \Rmn$ such that
	\begin{align*}
	X^\opts \in \frac{1}{1+\gamma}\svd_r(Z^\opts) \quad \textnormal{and} \quad \frac{X^\opts-Z^\opts}{\gamma} \in \partial g(X^\opts).
	\end{align*}
	Defining $R := \gamma^{-1}(Z^\opts - (1+\gamma) X^\opts)$ gives the equivalence for the non-convex Douglas-Rachford.
\end{proof}

\end{document}